\documentclass{article}
\usepackage{latexsym,amsmath,amssymb,amscd}
\begin{document}

\title{On the time evolution of a homogeneous isotropic gaseous universe }
\author{Tetu Makino \footnote{Professor Emeritus at Yamaguchi University, Japan /e-mail: makino@yamaguchi-u.ac.jp}}
\date{\today}
\maketitle

\begin{abstract}
We investigate the time evolution of a homogeneous isotropic universe
which consists of a perfect gas governed by the
Einstein-Euler-de Sitter equations.. Under suitable assumptions on the equation of state and assumptions on the present states, we give a mathematically rigorous proof
of the existence of the Big Bang at the finite past and expanding 
of the universe in the course to the
infinite future.
\end{abstract}

\newtheorem{Lemma}{Lemma}
\newtheorem{Proposition}{Proposition}
\newtheorem{Theorem}{Theorem}
\newtheorem{Definition}{Definition}
\newtheorem{Corollary}{Corollary}

\section{Introduction}

We consider the Einstein-Euler-de Sitter equations
\begin{align*}
& R_{\mu\nu}-\frac{1}{2}(g^{\alpha\beta}R_{\alpha\beta})g_{\mu\nu}-\Lambda g_{\mu\nu}=
\frac{8\pi G}{c^4}T_{\mu\nu}, \\
&T^{\mu\nu}=(c^2\rho +P)U^{\mu}U^{\mu}-Pg^{\mu\nu}
\end{align*}
for the metric
$$ds^2=g_{\mu\nu}dx^{\mu}dx^{\nu}. $$
Here the cosmological constant $\Lambda$ is supposed to be positive.
We assume\\

{\bf (A0) $P$ is a smooth function of $\rho>0$ such that
$0\leq P, 0\leq dP/d\rho <c^2$ for $\rho>0$ and $P \rightarrow 0$
as $\rho\rightarrow +0$.}\\

We consider the co-moving spherically symmetric metric
$$ds^2=e^{2F(t,r)}c^2dt^2-
e^{2H(t,r)}dr^2-
R(t,r)^2(d\theta^2+\sin^2\theta d\phi^2)$$
such that $U^0=e^{-F}, U^1=U^2=U^3=0$ for
$x^0=ct, x^1=r, x^2=\theta, x^3=\phi$.
Then the equations on the region where $\rho >0$ turn out to be

\begin{subequations}
\begin{eqnarray}
e^{-F}\frac{\partial R}{\partial t}&=&V \label{Aa} \\
e^{-F}\frac{\partial \rho}{\partial t}&=&-(\rho+P/c^2)
\Big(\frac{V'}{R'}+\frac{2V}{R}\Big) \label{Ab} \\
e^{-F}\frac{\partial V}{\partial t}&=&
-GR\Big(\frac{m}{R^3}+\frac{4\pi P}{c^2}\Big)+\frac{c^2\Lambda}{3}R + \nonumber \\
&&-\Big(1+\frac{V^2}{c^2}-
\frac{2Gm}{c^2R}-\frac{\Lambda}{3}R^2\Big)
\frac{P'}{R'(\rho+P/c^2)}
\label{Ac} \\
e^{-F}\frac{\partial m}{\partial t}&=&-\frac{4\pi}{c^2}R^2PV.
\label{Ad}
\end{eqnarray}
\end{subequations}
Here $X'$ stands for $\partial X/\partial r$. We put
\begin{equation}
m=4\pi\int_0^r\rho R^2R'dr \label{I}
\end{equation}
and the coefficients of the metric are given by
\begin{equation}
P'+(c^2\rho+P)F'=0 \label{III}
\end{equation}
and
\begin{equation}
e^{2H}=\Big(1+\frac{V^2}{c^2}
-\frac{2Gm}{c^2R}-\frac{\Lambda}{3}R^2\Big)^{-1}(R')^2.
\label{II}
\end{equation}
For derivation of these equations, see \cite{ssEE}, \cite{MisnerS}.\\

Now we are going to consider solutions of the form
$$R(t,r)=a(t)r,\qquad \rho=\rho(t), \qquad F(t,r)=0.$$
Of course we consider $a(t)>0$.
Then (\ref{Aa}) reads
\begin{equation}
V(t,r)=\frac{da}{dt}r, \label{1}
\end{equation}
(\ref{I}) reads
\begin{equation}
m(t,r)=\frac{4\pi}{3}a(t)^3r^3\rho(t), \label{2}
\end{equation}
(\ref{Ac}) reads
\begin{equation}
\frac{d^2a}{dt^2}=
\Big(-\frac{4\pi G}{3}(\rho+3P/c^2)+\frac{c^2\Lambda}{3}\Big)a, \label{3}
\end{equation}
(\ref{Ab}) reads
\begin{equation}
\frac{d\rho}{dt}=-3(\rho+P/c^2)\frac{1}{a}\frac{da}{dt}, \label{4}
\end{equation}
and (\ref{Ad}) is reduced to (\ref{4}). Therefore we have to consider only
(\ref{3}) and (\ref{4}).

Let us write (\ref{3})(\ref{4}) as the first order system
\begin{subequations}
\begin{eqnarray}
\frac{da}{dt}&=&\dot{a}, \label{FSa} \\
\frac{d\dot{a}}{dt}&=&\Big(-\frac{4\pi G}{3}(\rho+3P/c^2)+
\frac{c^2\Lambda}{3}\Big)a, \label{FSb} \\
\frac{d\rho}{dt}&=&-3(\rho+P/c^2)\frac{\dot{a}}{a}. \label{FSc}
\end{eqnarray}
\end{subequations}

We shall investigate the solution $(a(t), \dot{a}(t), \rho(t))$ of the
system (\ref{FSa})(\ref{FSb})(\ref{FSc}) which satisfies the initial
conditions
\begin{equation}
a(0)=a_0,\qquad \dot{a}(0)=\dot{a}_0, \qquad \rho(0)=\rho_0.
\end{equation}

Let $]t_-, t_+[, -\infty\leq t_-<0<t_+\leq+\infty,$ be the maximal interval of existence of the
solution in the domain
\begin{equation}
\mathcal{D}=\{ (a,\dot{a}, \rho)\  |\  0<a, 0<\rho \}.
\end{equation}\\

We shall use the variable $\rho^{\flat}$, a given function of $\rho$,
defined by
\begin{equation}
\rho^{\flat}=\exp\Big[\int^{\rho}\frac{d\rho}{\rho+P/c^2}\Big].
\end{equation}The assumption (A0) implies that $\rho\leq \rho+P/c^2\leq 2\rho$ so that $\rho^{\flat}\rightarrow+\infty\Leftrightarrow
\rho\rightarrow +\infty$, and $ \rho^{\flat}\rightarrow 0\Leftrightarrow
\rho\rightarrow 0$. The equation (\ref{FSc}) reads
$$\frac{d\rho^{\flat}}{dt}=-3\frac{\rho^{\flat}}{a}\frac{da}{dt},
$$
therefore it holds that
\begin{equation}
\frac{\rho^{\flat}}{\rho^{\flat}_0}=\Big(\frac{a_0}{a}\Big)^3, \label{MD2}
\end{equation}
where $\rho_0^{\flat}=\rho^{\flat}|_{\rho=\rho_0}$,
as long as the solution exists.\\

Let us recall the theory of A. Friedman \cite{Friedman} (1922). 
If (\ref{3})(\ref{4}) hold, then it follows that
the quantity
$$X:=\Big(\frac{da}{dt}\Big)^2-
\Big(\frac{8\pi G}{3}\rho+\frac{c^2\Lambda}{3}\Big)a^2$$
enjoys
$$\frac{dX}{dt}=0.$$
Therefore there should exist a constant $K$ such that
\begin{equation}
\Big(\frac{da}{dt}\Big)^2=\Big(\frac{8\pi G}{3}\rho+
\frac{c^2\Lambda}{3}\Big)a^2-c^2K. \label{FK}
\end{equation}
Of course the constant $K$ is determined by the initial condition as
$$K=\frac{1}{c^2}\Big(\Big(\frac{8\pi G}{3}\rho_0+\frac{c^2\Lambda}{3}\Big)(a_0)^2
-(\dot{a}_0)^2\Big).$$
By defining $K$ as this, we can write
$$e^{2H}=(1-Kr^2)^{-1}a^2$$
and the metric turns out to be
$$ds^2=c^2dt^2-a(t)^2\Big(
\frac{dr^2}{1-Kr^2}+r^2(d\theta^2+\sin^2\theta d\phi^2)\Big).
$$
The equation (\ref{FK}) is nothing but \cite[(6)]{Friedman} when
$\rho\propto a^{-3}$. \\ 

{\bf Remark 1.}  Note that if (\ref{4})(\ref{FK}) hold, then we can claim that (\ref{3}) holds as long as 
$da/dt\not=0$. But, when $da/dt=0$, or, $a(t)=\bar{a}=\mbox{Const.}>0$,
then $\rho(t)$ should be a constant, say, $\bar{\rho}$ and
\begin{equation}
4\pi G(\bar{\rho}+3\bar{P}/c^2)=c^2\Lambda \label{6}
\end{equation}
should hold by (\ref{3}). Then
$$e^{2H}=(1-k\bar{a}^2r^2)^{-1}\bar{a}^2$$
with $$k=\frac{4\pi G}{c^2}(\bar{\rho}+\bar{P}/c^2). $$
In this case, we can put $\bar{a}=1$ without loss of generality, and
the metric turns out to be
$$ds^2=c^2dt^2-\Big(
\frac{dr^2}{1-kr^2}+r^2(d\theta^2+
\sin^2\theta d\phi^2)\Big).
$$
Since $k>0$ provided that $\bar{\rho}>0$, this is the so-called `Einsteins's universe (1917)'.\\

The following rule is trivial:

\begin{Proposition}
If $(a,\dot{a},\rho)=(\phi_0(t), \phi_1(t), \psi(t))$ satisfies (\ref{FSa})(\ref{FSb})
(\ref{FSc}) for $t\in I$, $I$ being an open interval,
then
$(a,\dot{a},\rho)=(\phi_0(t+\Theta), \phi_1(t+\Theta), \psi(t+\Theta))$
satisfies (\ref{FSa})(\ref{FSb})(\ref{FSc}) for 
$t\in I-\Theta:=\{t | t+\Theta \in I \}$, 
$\Theta$ being a constant,
and 
$(a,\dot{a},\rho)=(\phi_0(-t), -\phi_1(-t), \psi(-t))$ satisfies
(\ref{FSa})(\ref{FSb})(\ref{FSc}) for $t \in -I:=\{t | -t \in I\}$.
\end{Proposition}

\section{Big Bang}

Hereafter we shall consider the solution such that 
\begin{equation}
\frac{da}{dt}\Big|_{t=0}=\dot{a}_0 >0, \label{P}
\end{equation}
which requires 
\begin{equation}
\Big(\frac{8\pi G}{3}\rho_0+\frac{c^2\Lambda}{3}\Big)(a_0)^2>c^2K.
\label{PP}
\end{equation}

In this section we investigate what will happen
when we prolong the solution to the left, or, to the past. We claim
\begin{Theorem}
Assume (A0). If the initial data satisfy (\ref{P}) and $\displaystyle \frac{d^2a}{dt^2}\Big|_{t=0}\leq 0$, that is, 
\begin{equation}
4\pi G(\rho_0+3P_0/c^2)\geq c^2\Lambda, \label{8}
\end{equation}
then $t_-$ is finite and $a(t)\rightarrow 0,
\rho(t)\rightarrow +\infty$ as $t\rightarrow t_-+0$.
\end{Theorem}

Proof. 
First we claim that $\dot{a}(t)>0$ for $t_-<\forall t \leq 0$. Otherwise,
there exists $t_1\in ]t_-, 0[$ such that
$\dot{a}(t)>0$ for $t_1<t\leq 0$ and $\dot{a}(t_1)=0$, since $\dot{a}(0)=\dot{a}_0>0$ is supposed ((\ref{P})).
Therefore
\begin{align*}
&\frac{d}{dt}(4\pi G(\rho+3P/c^2))=4\pi G\Big(1+\frac{3}{c^2}\frac{dP}{d\rho}\Big)\frac{d\rho}{dt}= \\
& =-4\pi G\Big(1+\frac{3}{c^2}\frac{dP}{d\rho}\Big)\cdot 3(\rho+P/c^2)
\frac{\dot{a}}{a}<0,
\end{align*}
for $t_1<\forall t\leq 0$,
and we see that (\ref{8}) implies
$$4\pi G(\rho+3P/c^2) \geq c^2\Lambda $$
for $t_1<\forall t\leq 0$. Therefore
$$\frac{d}{dt}\Big(\frac{8\pi G}{3}\rho+\frac{c^2\Lambda}{3}\Big)a^2=
-\frac{2}{3}(4\pi G(\rho+3P/c^2)-c^2\Lambda)a
\dot{a}\leq 0$$
and
$$
\Big(\frac{8\pi G}{3}\rho(t)+\frac{c^2\Lambda}{3}\Big)a(t)^2\geq
\Big(\frac{8\pi G}{3}\rho_0+\frac{c^2\Lambda}{3}\Big)a_0^2
$$
for $t_1<\forall t \leq 0$. We are supposing that 
$$\delta:=(\dot{a}_0)^2=\Big(\frac{8\pi G}{3}\rho_0+\frac{c^2\Lambda}{3}\Big)a_0^2-
c^2K >0.$$
Therefore,
$$
\Big(\frac{8\pi G}{3}\rho(t)+\frac{c^2\Lambda}{3}\Big)a(t)^2-
c^2K \geq \delta$$
for $t_1<\forall t \leq 0$, which implies
$\dot{a}(t)^2\geq \delta$
for $t_1<\forall t\leq 0$, a contradiction to $\dot{a}(t_1)=0$.
Hence we see $\dot{a}(t)^2>0$ for $t_-<\forall t \leq 0$.
Since $\dot{a}(0)>0$, we have $\dot{a}(t)>0$ for $t_-<\forall t\leq 0$ and,
repeating the above argument, we see 
$da/dt=\dot{a}(t) \geq \sqrt{\delta}$ for $t_-<\forall t \leq 0$. 
This implies $t_->-\infty$. Otherwise, if $t_-=-\infty$, we would have
$a(t)\leq a_0+\sqrt{\delta}t \rightarrow -\infty$ as
$t\rightarrow -\infty$, a contradiction. Since $da/dt=\dot{a}>0, d\rho/dt<0$, we have $a(t)\rightarrow a_-(\geq 0), \rho(t)\rightarrow
\rho_-(\leq +\infty)$ as $t\rightarrow t_-+0$.

We claim $\rho_-=+\infty$. Otherwise, if $\rho_-<+\infty$, we
have $a_->0$ thanks to (\ref{MD2}). Then the solution 
could be continued across
$t_-$ to the left. Hence it should be that $\rho_-=+\infty$ and
$a_-=0$ thanks to (\ref{MD2}).
This completes the proof. $\square$\\

Since the metric is
$$ds^2=c^2dt^2-
a(t)^2\Big(\frac{dr^2}{1-Kr^2}+r^2(d\theta^2+
\sin^2\theta d\phi^2)\Big),$$
we can consider the volume of the space at $t=\mbox{Const.}$ is
$a(t)^3\cdot 2\pi^2/K^{3/2}$, provided that $K>0$ and the space is the 3-dimensional
hypersphere of radius $a(t)/\sqrt{K}$. Since $a(t)\rightarrow 0$ and
$\rho(t)\rightarrow +\infty$ as $t\rightarrow t_-+0$, this is nothing but the so called `Big Bang'.\\

If $K\leq 0$, then the condition (\ref{8}) can be dropped, that is, we have
\begin{Corollary}
 Suppose (\ref{P}): $\dot{a}_0>0$. If $K\leq 0$, that is, if
\begin{equation}
\Big(\frac{8\pi G}{3}\rho_0+\frac{c^2\Lambda}{3}\Big)(a_0)^2
\leq (\dot{a}_0)^2, \label{Esc}
\end{equation}
then the conclusion of Theorem 1 holds even if (\ref{8}) is not satisfied.
\end{Corollary}
In fact, suppose (\ref{Esc}). Since we have
$\dot{a}(t)>0$ for $t_-<\forall t \leq 0$ a priori
thanks to $K\leq 0$ (see (\ref{FK})), we have
$$\frac{d\rho}{dt}=-3(\rho+P/c^2)\frac{\dot{a}}{a}<0$$
for $t_-<\forall t\leq 0$ a priori. Therefore $\rho(t)\rightarrow
\rho_-(\leq +\infty), a(t)\rightarrow a_-$ as $t\rightarrow t_-$. Suppose
$\rho_-<+\infty$. Note that $\rho_->\rho_0>0$. Then $a_->0$ thanks to (\ref{MD2}) and $\dot{a}(t)$ tends to a finite limit. Thus, if 
$t_->-\infty$, the solution could be continued to the left across $t_-$,
a contradiction. Hence $t_-=\infty$. Since $\rho_-<+\infty$ is supposed,
we should have $a_->0$. But (\ref{FK}) and $K\leq 0$ implies
$$\frac{da}{dt}\geq \sqrt{\frac{8\pi G\rho_0}{3}}a\geq\delta:=\sqrt{\frac{8\pi G\rho_0}{3}}a_->0,
$$
which implies $a\leq a_0+\delta t \rightarrow -\infty$ as
$t\rightarrow -\infty$, a contradiction. Therefore it should be the case
that $\rho_-=+\infty$. Then there exists $t_0\in ]t_-, 0[$ such that
$$ 4\pi G(\rho(t_0)+3P(t_0)/c^2) \geq c^2\Lambda,
$$
and the conclusion of Theorem 1 holds. $\square$\\

 More precise behavior of $a(t), \rho(t)$ as $t\rightarrow t_-+0$ can
be obtained by an additional assumption on the behavior of the function $\rho\mapsto P$ as $\rho \rightarrow +\infty$. Let us give an example.

The equation of state for neutron stars is given by
$$P=Ac^5\int_0^{\zeta}
\frac{q^4dq}{(1+q^2)^{1/2}},
\qquad
\rho=3Ac^3\int_0^{\zeta}(1+q^2)^{1/2}q^2dq.$$
See \cite[p. 188]{LandauL}.
Then we have
\begin{equation}
P=\frac{c^2}{3}\rho(1+[\rho^{-1/2}]_1) \label{9}
\end{equation}
as $\rho\rightarrow +\infty$. Here $[X]_1$ stands for a convergent power series of the form
$\sum_{k\geq 1}a_kX^k$.

Thus, generalizing this situation,  we assume\\

{\bf (A1) There are constants $\Gamma,\sigma$ such that
$1\leq\Gamma <2, 0<\sigma\leq 1$ and
\begin{equation}
P=(\Gamma-1)c^2\rho(1+O(\rho^{-\sigma}))
\end{equation}
as $\rho\rightarrow +\infty$.}\\

Then
integrating (\ref{FSc}) gives
$$a=a_1\rho^{-\frac{1}{3\Gamma}}(1+O(\rho^{-\sigma}))$$
as $\rho\rightarrow +\infty$ with a positive constant $a_1$ and
$$\rho=\rho_1a^{-3\Gamma}(1+O(a^{3\Gamma\sigma}))$$
as $a\rightarrow 0$ with $\rho_1=(a_1)^{3\Gamma}$. Then the equation
(\ref{FK}) turns out to be
$$\frac{da}{dt}=\sqrt{\frac{8\pi G\rho_1}{3}}
a^{\frac{-3\Gamma+2}{2}}(1+O(a^{\nu})),$$
where $\nu:=\min(3\Gamma\sigma, 3\Gamma-2)$. Solving this, we have
\begin{align}
a&=(6\pi\Gamma^2G\rho_1)^{\frac{1}{3\Gamma}}(t-t_-)^{\frac{2}{3\Gamma}}
(1+O((t-t_-)^{\frac{2\nu}{3\Gamma}})), \\
\rho&=\frac{1}{6\pi \Gamma^2G}(t-t_-)^{-2}(1+
O((t-t_-)^{\frac{2\nu}{3\Gamma}}))
\end{align}\\

{\bf Remark 2.}  Suppose that $K>0$. Then the space at $t=\mbox{Const.}$ has a finite volume and the total mass is
$$M(t)=\frac{4\pi}{3}a(t)^3K^{-3/2}\rho(t).$$
However, since (\ref{FSc}) implies 
$$\frac{d}{dt}a^3\rho =-\frac{3P}{c^2}a^2\frac{da}{dt}\not= 0,$$
provided that $P>0$, the total mass is not conserved during the time
evolution of our universe. Therefore when $P\not= 0$, instead, we should consider the variable
$\rho^{\flat}$.
Then the modified total mass
$$M^{\flat}(t)=\frac{4\pi}{3}a(t)^3K^{-3/2}\rho^{\flat}(t)$$
is constant with respect to $t$.\\

{\bf Remark 3.}  Note that the arguments of this section are still valid
even if $\Lambda=0$. If $\Lambda=0$, (\ref{8}) is satisfied for any $\rho_0>0$, and
we have Big Bang whenever (14): $\dot{a}_0>0$.

\section{Expanding universe}

Let us investigate the behavior of the solution when continued to the right (to the future) as long as possible. 
We claim

\begin{Theorem}
Assume (A0)(A1). If the initial data satisfy (\ref{P}) and 
$\displaystyle \frac{d^2a}{dt^2}\Big|_{t=0}\geq 0$, that is,
\begin{equation}
4\pi G(\rho_0+3P_0/c^2)\leq c^2\Lambda, \label{Ex}
\end{equation}
then $t_+=+\infty$ and $a(t)\rightarrow +\infty,
\rho(t)\rightarrow 0$ as $t\rightarrow +\infty$.
\end{Theorem}

Proof. 
First we claim that $\dot{a}(t)>0$ for $0\leq\forall t <t_+$. Otherwise,
there exists $t_1\in ]0, t_+[$ such that
$\dot{a}(t)>0$ for $0\leq \forall t<t_1$ and $\dot{a}(t_1)=0$, since $\dot{a}(0)=\dot{a}_0>0$ is supposed ((\ref{P})).
Therefore
\begin{align*}
&\frac{d}{dt}(4\pi G(\rho+3P/c^2))=4\pi G\Big(1+\frac{3}{c^2}\frac{dP}{d\rho}\Big)\frac{d\rho}{dt}= \\
& =-4\pi G\Big(1+\frac{3}{c^2}\frac{dP}{d\rho}\Big)\cdot 3(\rho+P/c^2)
\frac{\dot{a}}{a}<0,
\end{align*}
for $0\leq\forall t <t_1$,
and we see that (\ref{Ex}) implies
$$4\pi G(\rho+3P/c^2) \leq c^2\Lambda $$
for $0\leq\forall t<t_1$. Therefore
$$\frac{d}{dt}\Big(\frac{8\pi G}{3}\rho+\frac{c^2\Lambda}{3}\Big)a^2=
-\frac{2}{3}(4\pi G(\rho+3P/c^2)-c^2\Lambda)a
\dot{a}\geq 0$$
and
$$
\Big(\frac{8\pi G}{3}\rho(t)+\frac{c^2\Lambda}{3}\Big)a(t)^2\geq
\Big(\frac{8\pi G}{3}\rho_0+\frac{c^2\Lambda}{3}\Big)a_0^2
$$
for $0\leq\forall t <t_1$. We are supposing that 
$$\delta:=(\dot{a}_0)^2=\Big(\frac{8\pi G}{3}\rho_0+\frac{c^2\Lambda}{3}\Big)a_0^2-
c^2K >0.$$
Therefore,
$$
\Big(\frac{8\pi G}{3}\rho(t)+\frac{c^2\Lambda}{3}\Big)a(t)^2-
c^2K \geq \delta$$
for $0\leq\forall t <t_1$, which implies
$\dot{a}(t)^2\geq \delta$
for $0\leq\forall t<t_1$, a contradiction to $\dot{a}(t_1)=0$.
Hence we see $\dot{a}(t)^2>0$ for $0\leq \forall t<t_+$.
Since $\dot{a}(0)>0$, we have $\dot{a}(t)>0$ for $0\leq \forall t<t_+$ and,
repeating the above argument, we see 
$da/dt=\dot{a}(t) \geq \sqrt{\delta}$ for $0\leq\forall t<t_+$. 

We claim $t_+=+\infty$. In fact, suppose $t_+<+\infty$. Since 
$da/dt>0, d\rho/dt<0$, we have $a \rightarrow a_+(\leq +\infty),
\rho\rightarrow \rho_+(\geq 0)$ as $t\rightarrow t_+-0$.
Suppose $\rho_+>0$. Then (\ref{MD2}) implies that $a_+<+\infty$.
Then $\dot{a}(t)$ tends to a finite positive limit, and the solution could be continued to the left across
$t=t_+$, a contradiction. Therefore it should be the case that 
$\rho_+=0$, and (\ref{MD2}) implies that $a_+=+\infty$. Then, as $t\rightarrow t_+$, 
$$\Big(\frac{8\pi G}{3}\rho+\frac{c^2\Lambda}{3}\Big)a^2-c^2K
\sim \frac{c^2\Lambda}{3}a^2$$
and
$\displaystyle \frac{da}{dt}\leq Ca$ with a sufficiently large constant
$C$. Then it is impossible that $a$ blows up with finite $t_+$. Therefore
we know $t_+=+\infty$.

We see $a(t)\geq a_0+\sqrt{\delta}t \rightarrow +\infty$
as $t\rightarrow +\infty$. Then (\ref{MD2}) implies $\rho(t)\rightarrow 0$. 
This completes the proof. $\square$\\

In order to observe the asymptotic behaviors of
$a(t), \rho(t)$, we put the additional assumption:\\

{\bf (A2)  There is a constant $\gamma$ such that $1<\gamma$ and
$P=O(\rho^{\gamma})$ as $\rho\rightarrow 0$.}\\

Under this assumption we can fix the variable $\rho^{\flat}$ by
\begin{equation}
\rho^{\flat}=\rho\exp\Big[-\int_0^{\rho}
\frac{P/c^2\rho}{1+P/c^2\rho}\frac{d\rho}{\rho}\Big],
\end{equation}
and then we have 
$$\rho^{\flat}=\rho(1+O(\rho^{\gamma-1})$$
as $\rho\rightarrow 0$. Suppose the conditions of Theorem 2 hold. Then,
since
\begin{align*}
\rho&=\rho_1a^{-3}(1+O(a^{-3(\gamma-1)})) \\
&=O(a^{-3})
\end{align*}
as $a\rightarrow \infty$ with a positive constant $\rho_1$, we have
$$\Big(\frac{8\pi G}{3}\rho+\frac{c^2\Lambda}{3}\Big)a^2-c^2K=
\frac{c^2\Lambda}{3}a^2(1+O(a^{-2}))
$$
so that
$$\frac{da}{dt}=c\sqrt{\frac{\Lambda}{3}}a(1+O(a^{-2})).
$$
Hence
\begin{align}
a&=a_1e^{c\sqrt{\Lambda/3}t}(1+O(e^{-2c\sqrt{\Lambda/3}t})), \\
\rho&=\rho_1a_1^{-3}e^{-3c\sqrt{\Lambda/3}t}(
1+O(e^{-\nu c\sqrt{\Lambda/3}t}))
\end{align}
as $t\rightarrow +\infty$ with a suitable positive constant $a_1$. Here
$\nu=\min (2, 3(\gamma-1))$. \\

{\bf Remark 4.}  It follows from Theorem 1 and Theorem 2 that
the Einstein's universe (1917) :
$$a(t)=1,\qquad \dot{a}(t)=0,\qquad \rho(t)=\bar{\rho},
$$
with
$$4\pi G(\bar{\rho}+3\bar{P}/c^2)=c^2\Lambda
$$
is unstable. Actually the initial condition
$$a_0=1,\qquad \dot{a}_0>0, \qquad
\rho_0=\bar{\rho}
$$
satisfies the conditions of Theorem 2, even if $\dot{a}_0$ is arbitrarily
small. Then the solution is infinitely expanding universe when continued to the future. On the other hand, the initial condition
$$a_0=1, \qquad \dot{a}_0 <0, \qquad \rho_0=\bar{\rho}
$$
satisfies the conditions of Theorem 1, by converting the direction of the
time $t$, even if $|\dot{a}_0|$ is arbitrarily small. Then the solution
is a `Big Crunch' when continued to the future, that is, 
$t_+<+\infty$ and $a(t)\rightarrow 0, \rho(t)\rightarrow +\infty$
as $t\rightarrow t_+-0$. \\

If $K\leq 0$, then the condition (\ref{Ex}) of Theorem 2 can be dropped, that is, we have
\begin{Corollary}
 Suppose (\ref{P}): $\dot{a}_0>0$. If $K\leq 0$, that is, if (\ref{Esc}) holds,
then the conclusion of Theorem 2 holds even if (\ref{Ex}) is not satisfied
\end{Corollary}

In fact, suppose $K\leq 0$. Then we have $\dot{a}(t)>0, d\rho/dt <0$
a priori for $0\leq \forall t <t_+$. Thus
$\rho(t)\rightarrow \rho_+ (\geq 0), a(t)\rightarrow a_+ (\leq+\infty)$
as $t\rightarrow t_+$. Suppose $\rho_+>0$. Then $a_+<+\infty$
thanks to (\ref{MD2}). If $t_+<+\infty$, then the solution could be continued to the right beyond $t_+$, a contradiction. Hence $t_+=+\infty$.
Since $\rho_+>0$ is supposed, we should have $a_+<+\infty$. But 
(\ref{FK}) and $K\leq 0$ implies
$$\frac{da}{dt}\geq c\sqrt{\frac{\Lambda}{3}}a,$$
which implies $a\geq \mbox{Const.}e^{c\sqrt{\Lambda/3}t}\rightarrow
+\infty$ as $t\rightarrow +\infty$, a contradiction. Therefore we can claim that
$\rho_+=0$ and there exists $t_0\in ]0, t_+[$ such that
$$ 4\pi G(\rho(t_0)+3P(t_0)/c^2)\leq c^2\Lambda,
$$
and the conclusion of Theorem 2 holds. $\square$\\

Let us consider the case in which (\ref{Ex}) does not hold, that is,
$$4\pi G(\rho_0+3P_0/c^2) >c^2\Lambda. $$
Then we have

\begin{Corollary}
 Suppose (A0) and
$$4\pi G(\rho_0+3P_0/c^2)>c^2\Lambda. $$
Then there exists $\epsilon>0$
depending upon $a_0, \rho_0$ such that, if $0<\dot{a}_0\leq\epsilon$,
then $t_+<+\infty$ and $a(t)\rightarrow 0, \rho(t)\rightarrow +\infty$
as $t\rightarrow t_+-0$.
\end{Corollary}

Proof. Note that $$ \frac{d}{d\rho}(\rho+3P/c^2)=1+\frac{3}{c^2}\frac{dP}{d\rho}>0.$$
So we can take $\rho_*>0$ such that $\rho_*<\rho_0$ and
$$\delta:=4\pi G(\rho_*+3P_*/c^2)-c^2\Lambda >0  \qquad\mbox{with}\qquad
 P_*=P|_{\rho=\rho_*}. $$
Consider the solution
$(a^0(t), \dot{a}^0(t), \rho^0(t))$ with the initial data $(a_0, 0, \rho_0)$
and take $T>0$ sufficiently small so that the solution exists on $[0,T]$ 
and satisfies $a^0(t)<a_0\Theta^{-1/3}$ for $0\leq\forall t \leq T$,
where $\Theta$ is a positive number such that $\rho^{\flat}_*/\rho_0^{\flat}<\Theta <1$. 
Here we denote $\rho_0^{\flat}=\rho^{\flat}|_{\rho=\rho_0},
\rho_*^{\flat}=\rho^{\flat}|_{\rho=\rho_*}$. 
Then, taking a positive 
$\epsilon$ sufficiently small, we can suppose that the solution
$(a(t), \dot{a}(t), \rho(t))$ with the initial data $(a_0, \dot{a}_0, \rho_0)$
exists and satisfies $a(t)<a_0\Theta^{-1/3}$
for $0\leq\forall t\leq T$, provided that $0<\dot{a}_0\leq\epsilon$.
Moreover we can suppose $\displaystyle \epsilon<\frac{\delta}{3}a_0T$. We claim that 
there exists $t_0\in [0,T]$ such that $\dot{a}(t_0)<0$. Otherwise,
$\dot{a}(t)\geq 0$ and $a(t)\geq a_0$ for
$0\leq\forall t \leq T$. Since $a(t)<a_0\Theta^{-1/3}$ for 
$0\leq\forall t\leq T$, we have $\rho^{\flat}(t)>\rho_*^{\flat}$,
therefore $\rho(t)>\rho_*$ for
$0\leq\forall t\leq T$. Hence
$$\frac{d\dot{a}}{dt}=-\frac{1}{3}(4\pi G(\rho+3P/c^2)-c^2\Lambda)a <-\frac{\delta}{3}a$$
for $0\leq\forall t\leq T$, and
$$
\dot{a}(T)<\dot{a}_0-\frac{\delta}{3}\int_0^Tadt \leq \dot{a}_0-
\frac{\delta}{3}a_0T 
\leq \epsilon-\frac{\delta}{3}a_0T <0,
$$
a contradiction to
$\dot{a}(T)\geq 0$. Therefore there should exist $t_0\in [0,T]$ such that 
$\dot{a}(t_0)<0$, that is, the universe will become contracting. Moreover
since $\rho(t_0)>\rho_*$, we have
$$4\pi G(\rho(t_0)+3P(t_0)/c^2)>c^2\Lambda. $$
Then Theorem 1 can be applied by converting the direction of the time $t$. $\square$ \\

{\bf Remark 5.}  In the case of Corollary 3, there is $t_*$ such that
$0<t_*<t_+$, $\dot{a}(t)>0$ for $0\leq\forall t<t_*$ and $\dot{a}(t_*)=0$.
Then, by the uniqueness of the solution of (\ref{FSa})(\ref{FSb})(\ref{FSc}) we see that
$t_-=2t_*-t_+$ and
$a(t)=a(t_++t_--t)$ for $t_-<t<t_*$, that is, the solution $a(t)$ performs
a Big Bang at $t=t_-$ and a Big Crunch at $t=t_+$.\\

Let us consider the solution continued to the future supposing that $\dot{a}_0<0$. By converting the direction of $t$
we can apply Theorem 1 and Corollary 1, that is,
either if the condition (\ref{8}) :
$$4\pi G(\rho_0+3P_0/c^2)\geq\Lambda$$
holds or if $K\leq 0$:
$$\Big(\frac{8\pi G}{3}\rho_0+\frac{c^2\Lambda}{3}\Big)(a_0)^2\leq (\dot{a}_0)^2$$
holds, then we have a Big Crunch. When (\ref{8}) does not hold and $K>0$,
then we can use

\begin{Corollary}
Suppose (A0) and
$$4\pi G(\rho_0+3P_0/c^2)<c^2\Lambda. $$
Then there exists $\epsilon >0$ depending upon $a_0, \rho_0$ such that,
if $-\epsilon \leq \dot{a}_0<0$, then
$t_+=+\infty$ and $a(t)\rightarrow +\infty,
\rho(t)\rightarrow 0$ as $t\rightarrow +\infty$.
\end{Corollary}
 Proof. Take $\rho_*>0$ such that $\rho_0<\rho_*$ and
$$\delta:=c^2\Lambda-4\pi G(\rho_*+3P_*/c^2)>0.$$
It is possible. Consider the solution $(a^o(t), \dot{a}^o(t), \rho^o(t))$ with
the initial data $ (a_0, 0, \rho_0)$ and take $T>0$ such that
the solution exists and satisfies $a^o(t)>a_0\Theta^{1/3}$ for
$0\leq t\leq T$, where $\Theta$ is a positive number such that
$\rho_0^{\flat}/\rho_*^{\flat}<\Theta<1$.
Then, taking $\epsilon>0$ sufficiently small, we can suppose that
the solution $(a(t), \dot{a}(t), \rho(t))$ with the initial
data $(a_0, \dot{a}_0, \rho_0)$ exists and satisfies $a(t)>a_0\Theta^{1/3}$
for $0\leq t\leq T$, provided that
$-\epsilon \leq \dot{a}_0<0$. We can suppose that
$\displaystyle \epsilon < \frac{\delta}{3}a_0\Theta^{1/3}T$. 
We can claim that $\dot{a}(T)>0$. In fact, since
$a(t)>a_0\Theta^{1/3}$, we have
$\rho^{\flat}(t)<\rho_*^{\flat}$, therefore, $\rho(t)<\rho_*$ for
$0\leq t\leq T$. Hence
$$\frac{d\dot{a}}{dt}=-\frac{1}{3}(4\pi G(\rho+3P/c^2)-c^2\Lambda)a
>\frac{\delta}{3}a$$
and
$$\dot{a}(T)>\dot{a}_0+\frac{\delta}{3}\int_0^Tadt
\geq \dot{a}_0+\frac{\delta}{3}a_0\Theta^{1/3}T>0.$$
Now, since $\dot{a}(T)>0$ and
$$c^2\Lambda-4\pi G(\rho(T)+3P(T)/c^2)>0,$$
Theorem 2 can be applied to claim that
$t_+=+\infty$ and $a(t)\rightarrow +\infty, \rho(t)\rightarrow
0$ as $t\rightarrow +\infty$. $\square$\\

{\bf\Large Appendix}\\

Let us observe possible scenarios assuming that the equation of state is $P=0$. 

Integrating (\ref{FSc}) gives
$$\rho=\rho_1a^{-3}
$$
with a positive constant $\rho_1=\rho_0(a_0)^3$. The Friedman equation (\ref{FK}) turns out to be
$$\Big(\frac{da}{dt}\Big)^2=\frac{8\pi G\rho_1}{3}\frac{1}{a}+
\frac{c^2\Lambda}{3}a^2-c^2K.
$$
We shall consider the case $K>0$, and, suppose
$$0<\dot{a}_0<
\sqrt{\frac{8\pi G\rho_1}{3}\frac{1}{a_0}+\frac{c^2\Lambda}{3}(a_0)^2}.
$$
Let $]t_-, t_+[$ be the maximal interval of existence of the solution
in the domain$\{ 0<a, |\dot{a}|<\infty, 0<\rho \}$, and let
$]t_-^+, t_+^+[$ be that in the domain
$\{ 0<a, 0<\dot{a}, 0<\rho \}$.
When $t\in ]t_-^+, t_+^+[$, we can take
$$
\frac{da}{dt}=
\sqrt{\frac{8\pi G\rho_1}{3}\frac{1}{a}+\frac{c^2\Lambda}{3}a^2-c^2K},
$$
which can be integrated as
$$I:=\int \sqrt{ \frac{a}{\displaystyle a^3-\frac{3K}{\Lambda}a+\frac{8\pi G\rho_1}{c^\Lambda}} }da
=c\sqrt{\frac{\Lambda}{3}}(t-t_0),
$$
where $t_0$ is arbitrary. Introducing the variable $\xi$ and the 
parameter $\alpha$ defined by
$$a=\sqrt{\frac{K}{L}}\xi,\qquad
\alpha=\frac{4\pi G\rho_1}{c^2}\sqrt{\frac{\Lambda}{K^3}},
$$
we write
$$I=\int\sqrt{\frac{\xi}{\xi^3-3\xi +2\alpha}}d\xi.$$

Let us consider the following cases:

[\!( Case-0 )\!]: $0<\alpha<1$,

[\!( Case-1 )\!]: $\alpha=1$,

[\!( Case-2 )\!]: $1<\alpha <+\infty$. 

Let us denote
$$f_{\alpha}(\xi):=\xi^3-3\xi+2\alpha.$$
First we consider the case [\!( Case-1 )\!]: $\alpha =1$. In this case
we have
$$f_1(\xi)=(\xi-1)^2(\xi+2),$$
and
$$\sqrt{\frac{\xi}{\xi^3-3\xi+2}}=\frac{1}{|\xi-1|}\sqrt{\frac{\xi}{\xi+2}}.
$$
Therefore, using the variable $x$ defined by
$$x=\sqrt{\frac{\xi}{\xi+2}},$$
we see the integral $I$ is given as follows:

[\!( Case-1.0 )\!]: for $\displaystyle \frac{1}{\sqrt{3}}<x<1$
$$I=\frac{1}{\sqrt{3}}\log
\frac{\sqrt{3}x-1}{\sqrt{3}x+1}+
\log\frac{1+x}{1-x},$$
and

[\!( Case-1.1 )\!]: for $0<x<\displaystyle\frac{1}{\sqrt{3}}$
$$I=\frac{1}{\sqrt{3}}
\log\frac{1+\sqrt{3}x}{1-\sqrt{3}x}+
\log\frac{1-x}{1+x}.$$
The case [\!( Case-1.0 )\!] corresponds to $1<\xi$, and
the solution $a(t)$ exists and monotone increases on $]-\infty,
+\infty [$. Asymptotic behavior is:
\begin{align*}
&a(t) \sim Ce^{c\sqrt{\Lambda/3}t}\quad\mbox{as}\quad t\rightarrow +\infty, \\
&a(t)-\bar{a}\sim Ce^{c\sqrt{\Lambda}t}\quad\mbox{as}\quad
t\rightarrow -\infty
\end{align*}
with
$$\bar{a}=\Big(\frac{4\pi G\rho_1}{c^2\Lambda}\Big)^{1/3}.$$
Here and hereafter $C$ stands for various positive constants. \\

Note that $(a,\dot{a},\rho)=(\bar{a}, 0, \rho_0)$,
where $\rho_0=\displaystyle\frac{c^2\Lambda}{4\pi G}$, is the Einstein's static universe (1917).\\

On the other hand, [\!( Case-1.1 )\!] corresponds to
$0<\xi<1$, and the solution $a(t)$ exists and monotone
increases on $]t_-, +\infty[$, where $t_-$ is finite. Asymptotic behavior is:
\begin{align*}
&\bar{a}-a(t)\sim Ce^{-c\sqrt{\Lambda}t}\quad
\mbox{as}\quad t\rightarrow +\infty, \\
&a(t) \sim (6\pi G\rho_1)^{1/3}(t-t_-)^{2/3}
\quad\mbox{as}\quad t\rightarrow t_-+0.
\end{align*}

Next we consider the case [\!( Case-0 )\!]: $0<\alpha<1$. In this case
the polynomial $f_{\alpha}$ has three roots, say, $\xi_1, \xi_2, -(\xi_1+\xi_2)$, where $0<\xi_1<1<\xi_2$, so that
$$f_{\alpha}(\xi)=(\xi-\xi_1)(\xi-\xi_2)(\xi+\xi_1+\xi_2).$$ Therefore $\xi>0$ such that
$f_{\alpha}(\xi)>0$ are divided into two separated intervals:

[\!( Case-0.0 )\!]: $0<\xi<\xi_1$,

[\!( Case-0.1 )\!]: $\xi_2<\xi$.

\noindent The case [\!( Case-0.0 )\!] gives a solution $a(t)$
which exists and monotone increases on $]t_-, t_+^+[$,
where $t_-, t_+^+$ are finite. Put $t_*=t_+^+$. Asymptotic behavior is:
\begin{align*}
&a(t)=\sqrt{\frac{K}{\Lambda}}
\Big(\xi_1-
\frac{4c^2\Lambda(1-\xi_1^2)}{\xi_1}(t_*-t)^2+[(t_*-t)^2]_2\Big)
\quad\mbox{as}\quad t\rightarrow t_*-0, \\
&a(t)\sim (6\pi G\rho_1)^{1/3}
(t-t_-)^{2/3}\quad\mbox{as}\quad t\rightarrow t_-+0.
\end{align*}
So, we see that the solution $a(t)$ can be continued to the right
across $t=t_*$ as a monotone decreasing function as
$a(t)=a(2t_*-t)$
for $t>t_*$. Hence $t_+=2t_*-t_-$ and
$$a(t)\sim (6\pi G\rho_1)^{1/3}(t_+-t)^{2/3}\quad
\mbox{as}\quad t\rightarrow t_+-0.$$

In the same way we see that the case [\!( Case-0.1 )\!] gives a solution
$a(t)$ exists on $]-\infty, +\infty[$, monotone decreases on $]-\infty, t_*[$ and 
monotone increases on $]t_*, +\infty[$, where $t_*=t_-^+$ is
finite. Asymptotic behavior is
\begin{align*}
&a(t)\sim Ce^{c\sqrt{\Lambda/3}t}\quad\mbox{as}\quad t\rightarrow +\infty, \\
&a(t)\sim C'e^{-c\sqrt{\Lambda/3}t}\quad\mbox{as}\quad t\rightarrow -\infty,
\end{align*}
while
$$a(t)=\sqrt{\frac{K}{\Lambda}}
\Big(\xi_2+
\frac{4c^2\Lambda(\xi_2^2-1)}{\xi_2}(t-t_*)^2+
[(t-t_*)^2]_2\Big)
$$
as $t\rightarrow t_*$, and $a(t)=a(2t_*-t)$.

Finally we see for the case [\!( Case-2 )\!]: $1<\alpha$, the solution 
$a(t)$ exists and monotone increases on $]t_-, +\infty[$,
where $t_-$ is finite. Asymptotic behavior is:
\begin{align*}
&a(t)\sim Ce^{c\sqrt{\Lambda/3}t}\quad\mbox{as}\quad
t\rightarrow +\infty, \\
&a(t)\sim (6\pi G\rho_1)^{1/3}(t-t_-)^{2/3}
\quad\mbox{as}\quad t\rightarrow t_-+0.
\end{align*}
The scenario of this case is nothing but the so-called `Lema\^{i}tre model' of
\cite[pp.70-71, \S 5.1, (ii)]{Islam}, \cite[p.121, \S 10.8, Case 3(i)]{Bondi} 
(\cite{Lemaitre1927}, \cite{Lemaitre1931}).
In fact there exists a unique $t=t_m$ such that
$a(t)=\bar{a}$. Then $d\dot{a}/dt=0$ at $t=t_m$ and
$d\dot{a}/dt <0, ([ >0 ])$ for $t_-<t<t_m, ([ t_m<t<+\infty )]$ respectively.
Note that
$$\frac{d\dot{a}}{dt}\sim -C(t-t_-)^{-4/3} \quad
\mbox{as}\quad t\rightarrow t_-+0
$$
and
$$\frac{d\dot{a}}{dt}\sim C'e^{c\sqrt{\Lambda/3}t}
\quad\mbox{as}\quad t\rightarrow +\infty.
$$
Since
$$\dot{a}\sim C(t-t_-)^{-1/3}\quad\mbox{as}\quad t\rightarrow t_-+0$$
and
$$\dot{a}\sim C'e^{c\sqrt{\Lambda/3}t}\quad\mbox{as}\quad
t\rightarrow +\infty,
$$
we see that $\dot{a}(t)$ attains its positive minimum
at $t=t_m$. The time period around $t=t_m$ is called 
`coasting period' in the cosmological context.\\

Summing up, all possible scenarios with positive $K$ are:\\

AS $\nearrow$ EE, EC $\searrow$ AS, BB $\nearrow$ AS, AS $\searrow$ BC,
EC $\searrow\nearrow$ EE, BB $\nearrow\searrow$ BC, 

BB $\nearrow$ EE,
EC $\searrow$ BC \\

\noindent  except for the Einstein's static universe (1917). Here AS means `asymptotically steady state',
EC `exponentially contracting', BB `Big bang', EE `exponentially expanding', and BC means `Big Crunch'.

\end{document}